\documentclass[12pt]{article}
\usepackage{graphicx}
\usepackage{amssymb}
\usepackage{amsmath}
\usepackage{amsfonts}
\usepackage{amsthm}
\usepackage[english]{babel}
\usepackage[latin1,ansinew]{inputenc}
\usepackage{a4wide}
\usepackage{color}
\usepackage{float}
\usepackage{caption}
\usepackage{subcaption}
\newcommand{\be}{\begin{eqnarray}}
	\newcommand{\ee}{\end{eqnarray}}

\newtheorem{theo}{Theorem}

\newtheorem{lemma}{Lemma}

\newcommand{\R}{\mathbb R}

\newcommand{\eps}{\epsilon}

\title{Oscillating Shock Profiles in Relativistic Fluid Dynamics}
\author{ \it Valentin Pellhammer \thanks{Department of Mathematics and Statistics, University of Konstanz, 78457 Konstanz, Germany (valentin.pellhammer@uni-konstanz.de) }}

\begin{document}
	\maketitle

\begin{abstract}
	This note studies a model for relativistic pure-radiation fluids with viscosity that was recently proposed by Bemfica, Disconzi and Noronha, and shows that there are shock waves whose continuous shock profiles, if existing, are oscillating in any variables. This behavior differs significantly from the situation in classical fluid dynamics, in which canonical state variables are monotone along the shock profile. This paper shows that a generic local scenario established in \cite{P22} for a general context extends globally for this particular model.
\end{abstract}
	
	\section{Introduction}
	As a new model for the description of viscous relativistic pure-radiation fluids Bemfica, Disconzi and Noronha proposed in \cite{BDN} a system of four partial differential equations,
	\begin{equation}
	\partial_\beta(T^{\alpha\beta}+\Delta T^{\alpha\beta})=0 \label{PDE}
	\end{equation}
	for the four velocity $u^\gamma,\gamma=0,1,2,3$ and temperature $\theta$.
	The ideal part of the energy momentum tensor is given by 
	\[
	T^{\alpha\beta} = \frac{\partial(p(\theta)\psi^\alpha)}{\partial\psi_\beta}=\theta^3p'(\theta)\psi^\alpha\psi^\beta+p(\theta)g^{\alpha\beta},\quad p(\theta)=\frac{1}{3}\theta^4,
	\]
	written in Godunov variables $\psi^\gamma = \frac{u^\gamma}{\theta}$.
	The dissipative part is
	\[
	\Delta T^{\alpha\beta}=-B^{\alpha\beta\gamma\delta}(\psi)\frac{\partial \psi_\gamma}{\partial x^\delta},
	\]
	where 
	\[
	B^{\alpha\beta\gamma\delta}=\eta B_{visc}^{\alpha\beta\gamma\delta}-\mu B_{1}^{\alpha\beta\gamma\delta}-\nu B_{2}^{\alpha\beta\gamma\delta},
	\]
	with
	\[
	B_{visc}^{\alpha\beta\gamma\delta} = \Pi^{\alpha\gamma}\Pi^{\beta\delta} + \Pi^{\alpha\delta}\Pi^{\beta\gamma}-\frac{2}{3}\Pi^{\alpha\beta}\Pi^{\gamma\delta} ,\quad \Pi^{\alpha\beta}=u^{\alpha}u^\beta+g^{\alpha\beta},
	\]
	\[
	B_{1}^{\alpha\beta\gamma\delta} = (3u^\alpha u^\beta+\Pi^{\alpha\beta})(3u^\gamma u^\delta+\Pi^{\gamma\delta}),\quad B_{2}^{\alpha\beta\gamma\delta}=(u^\alpha{\Pi^\beta}_\eps+u^\beta{\Pi^\alpha}_\eps)(u^\gamma\Pi^{\delta\eps}+u^\delta\Pi^{\gamma\eps}).
	\]
	Additional to the viscosity parameter $\eta>0$ it depends on parameters $\mu, \nu > 0$ that were introduced to make the model causal. 
	The model \eqref{PDE} is a augmentation of the classical Eckart model (\cite{We}, Sec. 2.11), in order to overcome its lack of causality.
	In fact, it was shown in \cite{BDN,F21}, that all modes of the dissipation tensor are subluminal, if and only if
	\begin{equation}
	\mu\geq \frac{4}{3}\eta,\quad \nu\leq\left(\frac{1}{3\eta}-\frac{1}{9\mu}\right)^{-1}.\label{luminal}
	\end{equation}
	The non dissipative counterpart to \eqref{PDE} are given by the Euler equations
	\begin{align}\label{Euler}
		\partial_\beta T^{\alpha\beta}=0.
	\end{align}
	Fundamental solutions to the Euler equations are shock waves \cite{Lx}, i.e. discontinuous solutions of the prototypical form
	\begin{align}\label{shockWave}
	\psi_S(x)=\begin{cases}
		\psi_-,&x^\beta\xi_\beta<0\\\psi_+,&x^\beta\xi_\beta>0
	\end{cases},\quad \xi^\beta\xi_\beta=1,
	\end{align}
	and the question is whether these solutions correspond to solutions of the model with dissipation \eqref{PDE}. This is naturally done by means of \emph{dissipative shock profiles}.
	A shock wave \eqref{shockWave} has a \emph{dissipative shock profile} if there is a solution $\psi(x)=\tilde{\psi}(x^\beta\xi_\beta)$ to \eqref{PDE}
	that fulfills
	\begin{align}\label{endStates}
	\lim\limits_{z\to\pm\infty} \tilde{\psi}(z) = \psi_\pm.
	\end{align}
	Plugging $\tilde{\psi}$ into \eqref{PDE}, integrating once and using \eqref{endStates} yields that 
	the shock wave \eqref{shockWave} has a dissipative profile if and only if $\tilde{\psi}$ is a solution to 
	\begin{equation}
		\xi_\beta\xi_\delta B^{\alpha\beta\gamma\delta}(\psi)\psi'_\gamma = \xi_\beta T^{\alpha\beta}(\psi)-q^\alpha,\quad q^\alpha:=\xi_\beta T^{\alpha\beta}(\psi_\pm)\label{profileEquation}
	\end{equation}
	and \eqref{endStates}, i.e. \eqref{profileEquation} admits a heteroclinic profile connecting $\psi_-$ with $\psi_+$.
	
	It was shown in \cite{F21} that in the case where the second inequality on \eqref{luminal} holds strictly there are always Lax shocks that do not have a dissipative shock profile. Based on that it was suggested that one might prefer parameters with
	\begin{equation}
		\mu\geq \frac{4}{3}\eta,\quad \nu=\left(\frac{1}{3\eta}-\frac{1}{9\mu}\right)^{-1},\label{sharpCausal}
	\end{equation}
	in the model. The choice \eqref{sharpCausal} is referred to as \emph{sharply causal}, as in this scenario there exists a mode that is luminal. The purpose of this note is to provide further insight in the nature of shock profiles for this choice of parameters \eqref{sharpCausal}.
	
	In many known physical contexts shock profiles have been found to be always monotone in natural variables. This is distinctly not the case for the present model as the following result shows.
	
	\begin{theo}\label{Result2} For any values $\eta,\mu,\nu>0$ satisfying \eqref{sharpCausal},
		there is always a range of Lax-shocks that have either a oscillatory profile, i.e. it is not component wise monotone in any variables, or no profile at all.
	\end{theo}	
	A crucial point for the physical meaning of the shock profiles and even the model itself is the dynamical stability. Seen from this perspective the possible occurrence of oscillating shock profiles described in Theorem \ref{Result2} is rather delicate, as it was seen in other applications \cite{LLW,PSW}, that oscillatory behavior can indicate a lack of dynamical stability.
	
\section{The profile equations}
	\setcounter{equation}{0}
	
We follow \cite{F21} to simplify system \eqref{profileEquation}. We assume that \eqref{shockWave} is a standing shock wave with $(1,0,0)$ as the spacial direction of propagation of the shock front, i.e.
\[
(\xi_0,\xi_1,\xi_2,\xi_3)=(0,1,0,0).
\]
This choice causes no loss of generality by the Lorenz invariance and isotropy of the system. With this choice \eqref{profileEquation} reads 
\begin{align}\label{profileEquation2}
	B^{\alpha 1\gamma 1}(\psi)\psi_\gamma'=T^{\alpha 1} -q^\alpha.
\end{align}
We can now consider only states $\psi$ with $\psi^2=\psi^3=0$, which reduces \eqref{profileEquation2} to a planar dynamical system on the state space $\Psi=\{\psi=(\psi^0,\psi^1)\in\R^2:\psi^0>|\psi^1|\}$,
where the indices in \eqref{profileEquation2} run over $0$ and $1$ instead from $0$ to $3$.
The temperature and the velocity are now given as
\begin{align}\label{FormulasV}
\theta =(-\psi_\gamma\psi^\gamma)^{-\frac{1}{2}},\quad (u,v):=(u^0,u^1)=\theta(\psi^0,\psi^1)
\end{align}
with $u^2=1+v^2$.

\begin{lemma}\label{lemmaRestPoints}
	System \eqref{profileEquation2} has more than one rest point if and only if 
	\[
	q^1>0,\quad (q^1)^2<(q^0)^2<\frac{2}{\sqrt{3}}(q^1)^2.
	\]
	In this case it has exactly two rest points. These states form a standing  Lax-shock with respect to the Euler equations in right-moving or left-moving flow if $q^0>0$ or $q^0<0$ respectively. 
	\begin{proof}
		Cf. \cite{F21}, Lemma 1.
	\end{proof}
\end{lemma}

As system \eqref{odeSystem2} is invariant under the
reflection $q^0 \mapsto -q^0$, $u^1 \mapsto -u^1$ we can assume w.l.o.g.
\[
q^0>0,\quad q^1>0
\]
The pair of equilibria existing by Lemma \ref{lemmaRestPoints} depend on 
\[
\tilde{q}:=\left(\frac{q^1}{q^0}\right)^2 
\]
and is given by $(\psi_-,\psi_+)=(\psi_-(\tilde{q}),\psi_+(\tilde{q}))$ with
\begin{align}\label{VSquared}
\psi_\pm(\tilde{q}) = \left(\frac{4}{3}v^2+\frac{1}{3}\right)^\frac{1}{4}\begin{pmatrix}
	(1+v^2)^\frac{1}{2}\\v		
\end{pmatrix}\bigg\vert_{v=v_\pm(\tilde{q})},\quad 
v_\pm^2(\tilde{q}) = \frac{(2\tilde{q}-1)\mp \sqrt{\tilde{q}(4\tilde{q}-3)}}{4(1-\tilde{q})}.
\end{align}
They exist if and only if 
\begin{align*}
	\frac{3}{4}<\tilde{q}<1.
\end{align*}
For later use we note that 
\begin{align}\label{intervallsVsquare}
	v_+^2\in\left(\frac{1}{8},\frac{1}{2}\right).
\end{align}
We can furthermore assume w.l.o.g. $q^1 = 1$ and hence $\tilde{q}=(q^0)^{-2}$, cf. \cite{F21}.
Note that the limit $\tilde{q}\to \frac{3}{4}$ refers to the zero amplitude limit of the shock wave and $\tilde{q}\to 1$ refers to the infinite amplitude of the shock wave.

Consider parameters with \eqref{sharpCausal}.
A scaling of  the independent variable in \eqref{profileEquation2} according to $t\mapsto \frac{t}{\mu}$ and introducing the parameter
\begin{align}\label{defEps}
	\eps:= \frac{4\nu}{3\mu}
\end{align}
results in the two parameter dynamical system
\begin{align}\label{odeSystem2}
	B^\#(\psi,\eps)\psi' = F(\psi,\tilde{q})
\end{align}
with
\begin{align*}
	F(\psi,q) :=
	\begin{pmatrix}
		F^0\\F^1
	\end{pmatrix}
	:=
	\begin{pmatrix}
		-\frac{4}{3}\theta^4vu + q^0\\
		\theta^4(\frac{4}{3}v^2+\frac{1}{3})-q^1
	\end{pmatrix}
\end{align*}
and
\[
B^\#(\psi,\eps):=\eps B_{visc}(\psi) - B_{1}(\psi) - \frac{9\eps}{4-\eps} B_{2}(\psi),
\]	
where 
\[
B_{visc}(\psi) = \begin{pmatrix}
	u^2v^2 & -u^3v\\-u^3v & u^4
\end{pmatrix},\quad
B_{1}(\psi) =
\begin{pmatrix}
	16u^2v^2 & -4uv(4v^2+1)\\
	-4uv(4v^2+1) & (4v^2+1)^2
\end{pmatrix},
\]
\[
B_{2}(\psi) =
\begin{pmatrix}
	(u^2+v^2)^2 & -2(u^2+v^2)uv\\
	-2(u^2+v^2)uv & 4u^2v^2
\end{pmatrix}.
\]
With \eqref{defEps}, system \eqref{odeSystem2} fulfills the condition \eqref{sharpCausal} if and only if 
\begin{align*} \label{epsIntervall}
	0<\eps\leq 1.
\end{align*}
	\begin{lemma}\label{saddleAttactor}
	Consider system \eqref{odeSystem2} and assume $0<\eps\leq1$ as well as $\frac{3}{4}<\tilde{q}<1$. Then $\psi_-(\tilde{q})$ is a hyperbolic saddle and $\psi_+(\tilde{q})$ is an attractor.
	\begin{proof}
		Fix $0<\eps\leq1$ and	$\frac{3}{4}<\tilde{q}<1$.
		It was already shown in \cite{F21}, that $\psi_-(\tilde{q})$ is a saddle
		and $\det (DF(\psi_+(\tilde{q})))<0$
		for each choice of $\frac{3}{4}<\tilde{q}<1$.
		As
		\[
		\det(B^\sharp(\psi,\eps)) = \frac{9\eps((8+\eps)v^2+\eps-1)}{\eps-4},
		\]
		we find $\det(B^\#(\psi,\eps)) < 0$ whenever $v^2>\frac{1-\eps}{8+\eps}$.
		Looking at \eqref{intervallsVsquare} and observing $\frac{1-\eps}{8+\eps}\in[0,\frac{1}{8})$ for all $0<\eps\leq 1$ yields that $\det(B_\mu^\sharp(\psi_+(\tilde{q})))<0$ and therefore
		\[
		\det\left(B^\sharp(\psi_+(\tilde{q}),\eps)^{-1}\; DF(\psi_+(\tilde{q}))\right)>0.
		\] 
		This shows that $\psi_+(\tilde{q})$ is an attractor or a repeller. A simple calculation yields
		\[
		trace\left(B^\sharp(\psi,\eps)^{-1}DF(\psi)\right) = -\frac{3v}{(4-\eps)\det(B^\sharp(\psi,\eps))}((8+\eps^2)v^2+\eps^2-6\eps-4)
		\]
		which is negative if $0<v^2<\frac{-\eps^2+6\eps+4}{\eps^2+8}$. From the fact that
		$\frac{-\eps^2+6\eps+4}{\eps^2+8}\in(\frac{1}{2},1]$ for all $0<\eps\leq 1$ and \eqref{intervallsVsquare} we can deduce
		\[
		trace\left(B^\sharp(\psi_+(\tilde{q}),\eps)^{-1}DF(\psi_+(\tilde{q}))\right)<0,
		\]
		which shows that $\psi_+(\tilde{q})$ is an attractor.
	\end{proof}
\end{lemma}
We study \eqref{odeSystem2} as a dynamical system dependent on two parameters
in the parameter space
\[
\Omega:=\left\{(\eps,\tilde{q}):0<\eps\leq 1,\frac{3}{4}<\tilde{q}<1\right\},
\]
since it covers all possible shocks and choices of parameters with \eqref{sharpCausal} for the model.
	
	\section{Rest points with non real eigenvalues}\label{SecNotMonotone}
	\setcounter{equation}{0}
	Lemma \ref{saddleAttactor} makes sure that $\psi_+$ is an attractor
	for each choice of $\frac{3}{4}<\tilde{q}<1$. 
	As such it is either a stable focus or a stable node, depending on weather the eigenvalues of
	the lineraziation of \eqref{profileEquation2} have a nonzero imaginary part or not.
	A \emph{scaled} version of the linarization around a state $\bar{\psi}$ is given by
	\begin{align}\label{linearization}
		B^\#(\bar{\psi},\eps)\psi' = A(\bar{\psi})\psi
	\end{align}
	with 
	\[
	A(\psi) := \begin{pmatrix}
		v(6v^2+5) & -u(6v^2+1)\\ -u(6v^2+1) &3v(2v^2+1)
	\end{pmatrix}.
	\]
	
	Depending on the two parameters $(\eps,\tilde{q})\in\Omega$ the
	eigenvalues of the matrix $B^\#(\psi_+(\tilde{q}),\eps)^{-1}A(\psi_+(\tilde{q}))$ determine if $\psi_+(\tilde{q})$ is a 
	stable node or a stable spiraling focus. We therefore consider the discriminant of the characteristic polynomial of the matrix $B^\#(\psi,\eps)^{-1}A(\psi)$, which is given by
	\begin{align*}
	D(\psi,\eps) := trace(B^\#(\psi,\eps)^{-1}A(\psi))^2 - 4 \det(B^\#(\psi,\eps)^{-1}A(\psi))
	=\frac{1}{\det(B^\#(\psi,\eps))^2}\tilde{P}(\psi,\eps)
	\end{align*}
	with
	\[
	\tilde{P}(\psi,\eps) = trace(adj(B^\#(\psi,\eps))A(\psi))^2-4\det(B^\#(\psi,\eps))\det(A(\psi))),
	\]
	where $adj(B^\#)$ denotes the adjunct matrix of $B^\#$.
	A simple calculation yields
	\begin{align*}
		\det(B^\#(\psi,\eps)) = \frac{9\eps}{\eps-4}((\eps+8)v^2+\eps-1),\quad \det(A(\psi)) = 2v^2-1
	\end{align*}
	and
	\[
	trace(adj(B^\#(\psi,\eps))A(\psi)) = \frac{3v}{\eps-4}((8+\eps^2)v^2+\eps^2-6\eps-4),
	\]
	from which we deduce
	\[
	\tilde{P}(\psi,\eps) = \frac{9}{(\eps-4)^2}P(v^2,\eps)
	\]
	with the Polynomial $P$ given by
	\[
	P(z,\eps) = a_3(\eps)z^3+a_2(\eps)z^2 + a_1(\eps)z +a_0(\eps)
	\]
	where
	\begin{align*}
		a_0(\eps)=\eps(4\eps^2-20\eps+16),\quad a_1(\eps) = \eps^4 - 16\eps^3 + 84\eps^2 - 112\eps + 16\\
		a_2(\eps)= 2\eps^4 - 20\eps^3 - 24\eps^2 + 160\eps - 64,\quad a_3(\eps) = \eps^4 + 16\eps^2 + 64.
	\end{align*}
	For each $\eps>0$ and $\psi\in\Psi$, the sign of $D(\psi,\eps)$ coincides with the sign of $P(v^2,\eps)$ with $v$ given according to 
	\eqref{FormulasV}.

\begin{figure}[H]
		\begin{center}
			\includegraphics[trim =0mm 0mm 0mm 0mm ,scale=0.7]{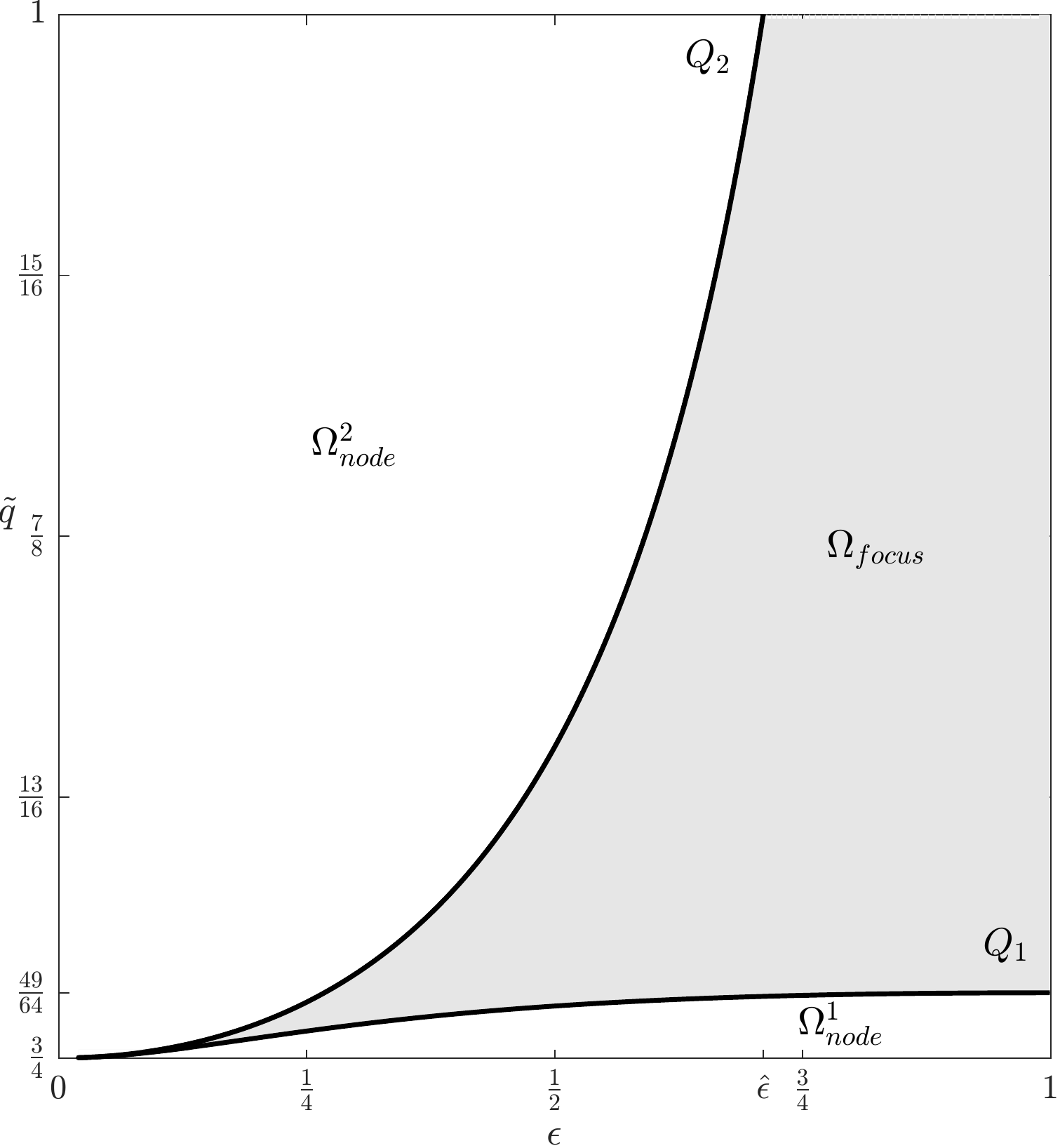}
			\caption{The parameter space for the equation \eqref{odeSystem2}. $\tilde{q}$ controls the strength of the shock and $\eps$ refers to the choice of dissipation. The rest point $\psi_+$ is a stable focus for $(\eps,\tilde{q})\in\Omega_{focus}$ and a stable node for
				$(\eps,\tilde{q})\in\Omega_{node}=\Omega_{node}^1\cup\Omega_{node}^2$. For every choice of sharply causal dissipation there are shock waves, whose profiles, if they exist, can only be oscillatory. \label{seperatrices}}
		\end{center}
	\end{figure}

	\begin{lemma}\label{lemmaNst}
		Considers $P$ as a parameter dependent polynomial in $z$ with parameter $0\leq\eps\leq 1$. 
		For each $0\leq\eps\leq 1$, $P$ has three real zeros $w_1(\eps)<w_2(\eps)\leq w_3(\eps)$ with $w_2(\eps)= w_3(\eps)$ if and only if $\eps=0$ and the following properties.
			\begin{enumerate}
	\item[(i)] $w_1(\eps)<0$ for all $0<\eps\leq 1$ and $w_1(0)=0$.
	\item[(ii)] for  $\hat{\eps}:= \frac{2}{3}\left(3\sqrt{6}-2\sqrt{16-6\sqrt{6}}-4\right)\approx 0.7103$ one has $ w_2((0,\hat{\eps}))\subset (\frac{1}{8},\frac{1}{2})$ and $w_2((\hat{\eps},1))\subset (0,\frac{1}{8})$ as well as $w_2(0)=\frac{1}{2}$,$w_2(\hat{\eps})=\frac{1}{8}$ and $w_2(1)=0$.
	\item[(iii)] $w_3(\eps)\in(\frac{1}{3},\frac{1}{2})$ for $0<\eps<1$ and $w_3(0)=\frac{1}{2}$, $w_3(1)=\frac{1}{3}$.

\end{enumerate}
\begin{proof}
	The discriminant of $P$ as a polynomial in $z$ is given by $D(\eps) = 1296\eps^5(4-\eps)^3\tilde{D}(\eps)$
	where
	\[
	\tilde{D}(\eps) = -4\eps^5 + 179\eps^4 - 844\eps^3 + 880\eps^2 - 32\eps + 64.
	\]
	Setting $\zeta=\frac{1}{\eps}$ yields
	\[
	\tilde{D}(\eps) =\zeta^{-5}\tilde{D}_2(\zeta),\quad \text{with}  \quad\tilde{D}_2(\zeta)= 64\zeta^5 - 32\zeta^4 + 880\zeta^3 - 844\zeta^2 + 179\zeta - 4.
	\]
	The polynomial
	\[
	\tilde{D}_2(x+1) = 64x^5 + 288x^4 + 1392x^3 + 2244x^2 + 1323x + 243
	\]
	has only positive coefficients and can thus not have a non negative real root. Therefore $\tilde{D}_2$ has no root in $[1,\infty)$ and has a constant sign in this interval. Since the leading coefficient of $\tilde{D}_2$ is positive one has $\tilde{D}_2(\zeta)>0$ for all $1<\zeta<\infty$ and we can conclude that  $D(\eps)>0$ for all $0<\eps<1$.
	This shows the existence of three distinct real roots $w_1(\eps)<w_2(\eps)<w_3(\eps)$, which depend continuously on $\eps$.
	
	One has
	\[
	P(z,1) = 27z(3z^2+2z-1),
	\]
	and therefore $w_1(1)=-1,w_2(1)=0$ and $w_3(1)=1/3$. 
	
	Since the leading coefficient of $P$ is positive and $P(0,\eps)=4\eps(\eps^2 - 5\eps + 4)>0$ for all $0<\eps\leq 1$
	it follows that $P$ has a negative real root, i.e. $w_1(\eps)<0$ for all $0<\eps\leq 1$. This shows assertion (i).

	As $P$ is cubic in $z$ with a positive leading coefficient and $P(\frac{1}{2},\eps)=\frac{9}{8}\eps^2(\eps - 4)^2>0$ we know that $w_3(\eps)\neq \frac{1}{2}$ for all $0<\eps\leq 1$. Since additionally $w_3(1)=\frac{1}{3}<\frac{1}{2}$ we can deduce that $w_3(\eps)<\frac{1}{2}$ for all $0<\eps\leq 1$ by the continuity of $w_3$.
	
	One has $P(\frac{1}{3},\eps)=\frac{16}{27}(\eps - 1)^2(\eps^2 - 4\eps + 1)$ and thus $P(\frac{1}{3},\eps)=0$ for a $0<\eps\leq 1$ if and only if $\eps\in\{2-\sqrt{3},1\}$ and it follows that
	$w_2(2-\sqrt{3})=\frac{1}{3}<w_3(2-\sqrt{3})$.
	Since additionally $w_3(\eps)=\frac{1}{3}$ if and only if $\eps=1$ we follow by continuity of $w_3$ that $\frac{1}{3}<w_3(\eps)$ for all $0<\eps<1$ which thus shows (ii).
	
	For $\eps>0$ one has
	\begin{align*}
	P\left(\frac{1}{8},\eps\right) &= \frac{9}{512}(9x^4+96x^3-560x^2+256x+64) \\ &=\frac{9\eps^2}{512}\left(9\eps^2+96\eps^1-560+256\frac{1}{\eps}+64\frac{1}{\eps^2}\right)\\
	&=\frac{9\eps^2}{512}\left(9\left(\eps+\frac{8}{3\eps}\right)^2+96 \left(\eps+\frac{8}{3\eps}\right)-608\right)
	\end{align*}
	and hence $P(\frac{1}{8},\eps)=0$ can easily be solved for $y:=(\eps+\frac{8}{3\eps})$. Doing this yields that $P(\frac{1}{8},\eps)=0$ for a $0<\eps\leq1$, i.e. $w_2(\eps)=\frac{1}{8}$  if and only if $\eps = \hat{\eps}:= \frac{2}{3}\left(3\sqrt{6}-2\sqrt{16-6\sqrt{6}}-4\right)$. Since $w_2(1)=0$ the continuity of  $w_2$ yields $w_2(\eps)<\frac{1}{8}$ for each $\hat{\eps}<\eps\leq 1$.
\end{proof}
	\end{lemma}
	The following Theorem reveals the transition of $\psi_+$ from a stable node to a stable focus for varying parameters in detail. In the case where $\psi_+$ is a focus a possibly existing orbit connecting $\psi_-$ with $\psi_+$ has a oscillating character due to its spiraling around $\psi_+$ at that end state, which implies Theorem \ref{Result2}.
	\begin{theo}
		The parameter range $\Omega$ decomposes into two open subsets $\Omega_{node}$ and $\Omega_{focus}$ and two seperatrices $Q_1$ and $Q_2$ such that the equilibrium $\psi_+$ is a node if and only if $(\eps,\tilde{q})\in\Omega_{node}$ and $\psi_+$ is a stable focus if and only if $(\eps,\tilde{q})\in\Omega_{focus}$. The sepreatrices $Q_1$ and $Q_2$ are the graphs of two functions $\tilde{q}_1$ and $\tilde{q}_2$ given by
		\[
		\tilde{q}_1:(0,1],\mapsto(\tfrac{3}{4},1),\quad \tilde{q}_2:(0,\hat{\eps}),\mapsto(\tfrac{3}{4},1).
		\]	
		 $\Omega_{focus}$ is the open subset contained between the two curves and $\Omega_{node}$ consists of the two open sets $\Omega_{node}^1$ below $Q_1$ and $\Omega_{node}^2$ above $Q_2$, cf. Figure \ref{seperatrices}.
		 Furthermore
		 \begin{align}\label{limits}
		 \tilde{q}_1(1)=\frac{49}{64},\quad\lim\limits_{\eps\to \hat{\eps}}\tilde{q}_2(\eps)=1,\quad \lim\limits_{\eps\to 0}\tilde{q}_1(\eps)=\lim\limits_{\eps\to 0}\tilde{q}_2(\eps)=\frac{3}{4}
		 \end{align}
		 \begin{proof}
		 	Consider the discriminant of the characteristic polynomial of the linearization at $\psi_+(\tilde{q})$ given by
		 	\[
		 	D(\psi_+(\tilde{q}),\eps)=\frac{9}{\det(B)^2(\eps-4)^2}P(v_+^2(\tilde{q}),\eps).
		 	\] 		 	
		 	and the function
		 	\[
		 	V_+:\left(\tfrac{3}{4},1\right)\to \left(\tfrac{1}{8},\tfrac{1}{2}\right),\quad\tilde{q}\mapsto v_+^2(\tilde{q})
		 	\]
		 	that links the value $\tilde{q}$ to the value for $v^2_+(\tilde{q})$ of the equilibrium $\psi_+(\tilde{q})$ via \eqref{VSquared}.
		 	The function $V_+$ is onto and strictly decreasing and has thus a continuous inverse $V_+^{-1}$ that is also strictly decreasing. 	 	
		 	We define
		 	\begin{align*}
		 	\tilde{q}_{1}(\eps):=(V_+^{-1}\circ w_{3})(\eps),\eps\in(0,1],\quad
		 	\tilde{q}_{2}(\eps):=(V_+^{-1}\circ w_{2})(\eps), \eps\in(0,\hat{\eps})
		 	\end{align*}
		 	where $w_2$ and $w_3$ are the functions defined in Lemma \ref{lemmaNst}. Since $V_+^{-1}$ is monotone one has $\tilde{q}_1<\tilde{q}_2$ on $(0,\hat{\eps})$ and the assertions \eqref{limits} easily follow using Lemma \ref{lemmaNst}.
	 	
		 	Since $P$ is cubic and the leading coefficient of $P$ is positive, Lemma \ref{lemmaNst} shows that $P(z,\eps)<0$ for $z\in(\frac{1}{8},\frac{1}{2})$ and $\eps\in(0,1]$ if and only if either $\eps\in(0,\hat{\eps})$ and $z\in(w_2(\tilde{q}),w_3(\tilde{q}))$ or $\eps\in[\hat{\eps},1)$ and $z>w_2(\tilde{q})$.
		 	
		 	For $\eps\in(0,\hat{\eps})$ and $q_1(\eps)<\tilde{q}<q_2(\eps)$ one has $v_+^2(\tilde{q})\in(w_2(\tilde{q}),w_3(\tilde{q}))$ and therefore $D(\psi_+(\tilde{q}),\eps)<0$ such that the eigenvalues of the linearization \eqref{linearization} have a nonzero real part and $\psi_+(\tilde{q})$ is a stable focus. If $\eps\in(0,\hat{\eps})$ and either $\tilde{q}<q_1(\eps)$ or $q_2(\eps)<\tilde{q}$ one has
		 	$v_+^2(\tilde{q})\notin(w_2(\tilde{q}),w_3(\tilde{q}))$ and thus $D(\psi_+(\tilde{q}),\eps)>0$. In this case the eigenvalues of \eqref{linearization} have a non negative imaginary part and $\psi_+(\tilde{q})$ is a stable node. 
		 	
		 	In the case $\eps\in(\hat{\eps},1]$ the same arguments yield that $\psi_+(\tilde{q})$ is a focus if $q_1(\eps)<\tilde{q}$ and a node if $q_1(\eps)<\tilde{q}$.
		 \end{proof}
	\end{theo}
	
\end{document}